\documentclass[12pt]{article}
\usepackage{amssymb}

\newcommand{\R}{\mathbb{R}}

\newtheorem{definition}{Definition}
\newtheorem{theorem}{Theorem}
\newtheorem{remark}{Remark}
\newtheorem{corollary}{Corollary}
\newtheorem{proposition}{Proposition}

\begin{document}

\title{Permanence criteria for Kolmogorov systems with delays}
\date{}
\author{Zhanyuan Hou}

\maketitle

{\centerline{School of Computing, London Metropolitan University,}
\centerline{166-220 Holloway Road, London N7~8DB, UK}
\centerline{Email: z.hou@londonmet.ac.uk}}

\begin{abstract}
In this paper, a class of Kolmogorov systems with delays are
studied. Sufficient conditions are provided for a system to have a compact uniform attractor. Then Jansen's result (J. Math. Biol. \textbf{25} (1987) 411--422) for autonomous replicator and Lotka-Volterra systems has been extended to delayed nonautonomous Kolmogorov systems with periodic or autonomous Lotka-Volterra subsystems. Thus, simple algebraic conditions are obtained for partial permanence and permanence. An outstanding feature of all these results is that the conditions are irrelevant of the size and distribution of the delays.
\end{abstract}



\textbf{Note} This paper is to be published in Proceedings of Royal Society of Edinburgh.

\section{Introduction}

Kolmogorov systems of differential equations have been used to model many biological problems and many variations of such systems have been extensively studied. Among the various investigations, permanence or uniform persistence is related to the problem of coexistence of species and received much attention in the last few decades. Here we only mention a few examples. Mierczy\'{n}ski and Schreiber \cite{MiSc} dealt with autonomous Kolmogorov system with robustly permanent subsystems. Kuang \cite{KuTa} and Tang \cite{TaKu} investigated delayed nonautonomous Kolmogorov systems and obtained permanence criteria, which depend on successful construction of Lyapunov functions or functionals. Yang \cite{Ya} studied persistence of a single-species Kolmogorov equation with delays. Examples of permanence of special classes of Kolmogorov systems without delays, including  Lotka-Volterra differential systems, are given by Ahmad and Lazer \cite {AhLa4}, Ahmad and Stamova \cite{AhSt}, Baigent and Hou \cite{BaHo}, Hofbauer and Schreiber \cite{HoSc}, Zhao and Jiang \cite{ZhJi}. In particular, for the autonomous Lotka-Volterra system
\begin{equation}\label{e1.1}
x'_i = x_i(r_i + A_ix), \quad i\in I_N=\{1, 2, \ldots, N\},
\end{equation}
where $A_i$ is the $i$th row of the $N\times N$ matrix $A$, Jansen \cite{Ja} (see also \cite[Ch.13]{HoSi}) proved that (\ref{e1.1}) is permanent if there is a vector $q\in \textrm{int}{\R}^N_+$ such that the inequality $q^T(r+A\hat{x})>0$ holds for every fixed point $\hat{x}\in\partial{\R}^N_+$. Examples of permanence for special delayed Kolmogorov systems are given by Chen, Lu and Wang \cite{ChLW}, Hou \cite{Ho3}--\cite{Ho5}, Li and Teng \cite{LiTe}, Liu and Chen \cite{LiCh}, Lu, Lu and Enatsu \cite{LuLuEn}, Mukherjee \cite{Mu}, Teng \cite{Te}, and the references therein. In particular, for autonomous Lotka-Volterra differential systems with multiple delays, sufficient conditions for permanence, which are easily checkable inequalities involving the constant coefficients of the system, were obtained in \cite{LuLuEn}.

In this paper, we are going to extend Jansen's result to a class of delayed nonautonomous Kolmogorov systems with Lotka-Volterra subsystems on the boundary $\partial{\R}^N_+$ having constant interactions and periodic intrinsic growth rates.

We shall consider the delayed nonautonomous Kolmogorov system
\begin{equation}
x'_i = x_if_i(t, x_t), \quad i\in I_N, \label{e1.2}
\end{equation}
where $f: {\R}_0\times C^+\to {\R}^N$ is continuous, ${\R}_0 = {\R}$ or $[\bar{t}, \infty)$ or $(\bar{t}, \infty)$ for some $\bar{t}\in {\R}$, $C^+ = C([-\tau, 0], {\R}^N_+)$ for some $\tau >0$ and $x_t(\theta)=x(t+\theta)$ for $\theta\in[-\tau, 0]$. Assume that the $f_i(t, \varphi)$ are locally Lipschitzian in $\varphi$. Then, for each $t_0\in{\R}_0$ and every $\varphi\in C^+$, the solution $x(t, t_0, \varphi)$ of (\ref{e1.2}) with $x_{t_0} = \varphi$ is unique, satisfies $x_t\in C^+$ for $t$ in its existing interval $[t_0, \mu)$, and depends on the initial data $(t_0, \varphi)$ continuously. In general, we may not have $\mu = +\infty$. Even if we have $\mu=+\infty$, $x(t, t_0, \varphi)$ may not be bounded on $[t_0, +\infty)$. Even if each solution exists on $[t_0, +\infty)$ and is bounded, (\ref{e1.2}) may not be uniformly bounded. A system is called \textit{uniformly bounded} if there is an $M>0$ such that every solution satisfies $|x(t)|<M$ for sufficiently large $t$. Moreover, the delays in (\ref{e1.2}) may cause dramatic changes on the behaviour of solutions. For example \cite{ChLW}, replacement of a term $x(t)$ by $x(t-\tau)$ in a uniformly bounded system may induce unbounded solutions. A system is called \textit{permanent} if there are $\delta>0$ and $M>\delta$ such that every solution in $\textrm{int}{\R}^N_+$ satisfies
\[
\forall i\in I_N, \; \forall\; \textup{large enough}\; t,\; \delta < x_i(t) < M.
\]
We say that solutions of (\ref{e1.2}) in $\textrm{int}{\R}^N_+$ are uniformly bounded away from the boundary $\partial{\R}^N_+$ if there is a $\delta>0$ such that
\[
\forall i\in I_N, \; \forall\; \textup{large enough}\; t,\; x_i(t) >\delta.
\]
Thus, (\ref{e1.2}) is permanent if and only if (\ref{e1.2}) is uniformly bounded and the solutions in $\textrm{int}{\R}^N_+$ are uniformly bounded away from $\partial{\R}^N_+$. Some available results (such as \cite{KuTa}, \cite{TaKu}, \cite{Mu} and \cite{MiSc}) on permanence assume the uniform boundedness while others (such as \cite{Te} and \cite{Ya}) prove the uniform boundedness under certain assumptions. There are various conditions for uniform boundedness of autonomous Lotka-Volterra systems (see \cite[Ch.15]{HoSi}).

In this paper, we are going to provide sufficient conditions to ensure that (\ref{e1.2}) has a \textit{compact uniform attractor}, an even better property than uniform boundedness, to be defined later. We shall see that our permanence results reply on this property. So we shall establish a few results for this property and then deal with permanence.

In this paper, the norm $|\cdot|$ on ${\R}^N$ is taken to be $|x| = \max\{|x_i|: i\in I_N\}$ and the norm $\|\cdot\|$ on $C^+$ is taken to be $\|\varphi\| = \max\{|\varphi(\theta)|: \theta\in [-\tau, 0]\}$.

\section{Main results}

We first describe a compact global attractor of (\ref{e1.2}) in $C^+$ that is positively invariant and every solution $x_t(t_0, \varphi)$ will enter and stay in this attractor after a finite time uniformly for $t_0\in {\R}_0$. Also, $\varphi$ in the boundary of this attractor if and only if $\varphi(\theta)\in\partial{\R}^N_+$ for all $\theta\in [-\tau, 0]$.
\begin{definition}\label{Def1} \textup{
System (\ref{e1.2}) is said to have a compact uniform attractor $\Omega\subset C^+$ if
\begin{itemize}
\item[(i)] $\Omega$ is compact;
\item[(ii)] for each $\varphi\in C^+$, there is a $T(\varphi)>0$ such that $x_t(t_0, \varphi)\in \Omega$ holds for all $t_0\in{\R}_0$ and all $t\geq t_0+T(\varphi)$;
\item[(iii)] $x_t(t_0, \varphi)\in\Omega$ holds for all $t_0\in{\R}_0$, $\varphi \in \Omega$ and $t\geq t_0$;
\item[(iv)] for each $\varphi\in\Omega$ and every $i\in I_N$, $\varphi_i(\theta_0)=0$ for some $\theta_0\in [-\tau, 0]$ if and only if $\varphi_i(\theta)\equiv 0$ on $[-\tau, 0]$.
\end{itemize}}
\end{definition}

Clearly, if (\ref{e1.2}) has a compact uniform attractor then it is uniformly bounded. Conversely, If (\ref{e1.2}) if uniformly bounded, does it have a compact uniform attractor? We cannot answer this question in general but the answer is positive if (\ref{e1.2}) is autonomous.

\begin{proposition}\label{Pro1}
If (\ref{e1.2}) is autonomous and $f$ is bounded on any bounded set $S\subset C^+$, then it is uniformly bounded if and only if it has a compact uniform attractor.
\end{proposition}
\textbf{Proof}
Suppose (\ref{e1.2}) is autonomous and uniformly bounded. Then
\begin{equation}\label{e2.0}
\exists M>0, \forall \varphi\in C^+, \exists T=T(\varphi)>0 \;\textup{such that}\; \forall t\geq T, |x(t, \varphi)|<M.
\end{equation}
Thus, for each $\varphi\in C^+$ with $\|\varphi\|\leq M$, there is $t_1(\varphi)\in [0, T(\varphi)+\tau)$ such that $\|x_t(\varphi)\|<M$ for all $t>t_1$ but, if $t_1>0$, $\|x_{t_1}(\varphi)\|=M$. Let
\[
S = \{x_t(\varphi): \varphi\in C^+, \|\varphi\|\leq M, t\geq t_1(\varphi)+\tau\}.
\]
We show that $\Omega = \bar{S}$, the closure of $S$, is a compact uniform attractor. Since $\Omega$ is bounded, by the assumption on $f$ there is a $\rho>0$ such that $|f(\varphi)|\leq \rho$ for all $\varphi\in \Omega$. Then, for $t\geq t_1(\varphi)$, we have
\[
\forall i\in I_N, x_i(t, \varphi) = x(t_1, \varphi)\exp\biggl(\int^t_{t_1}f_i(s, x_s(\varphi))ds\biggr)
\]
so
\[
\forall i\in I_N, x(t_1, \varphi)e^{-\rho(t-t_1)}\leq x_i(t, \varphi) \leq x(t_1, \varphi)e^{\rho(t-t_1)}.
\]
This shows that the functions over $[-\tau, 0]$ in $S$ are equicontinuous. As the inequalities for boundedness and equicontinuity of the functions in $S$ are retained for any limit function of a convergent sequence (refer to the proof of Theorem \ref{The2.1} given in section 3), by Arzela-Ascoli theorem, $\Omega$ is relatively compact. As $\Omega$ is also closed, $\Omega$ is compact. By the definition of $\Omega$, it is positively invariant. From (\ref{e2.0}) we see that for each $\varphi\in C^+$, $x_t(\varphi)\in \Omega$ for all $t\geq T(\varphi)+2\tau$. From the definition of $\Omega$ again, we know that for any $\psi\in \Omega$, $\psi_i(\theta_0)=0$ for some $\theta_0\in [-\tau, 0]$ if and only if $\psi_i(\theta)\equiv 0$. Therefore, $\Omega$ is a compact uniform attractor.
\rule{2mm}{3mm}

\noindent \textbf{Open Problem} If $f$ is bounded on ${\R}_0\times S$ for any bounded set $S\subset C^+$ and $f(t, \varphi)$ is (i) periodic in $t$ or (ii) almost periodic in $t$ or (iii) asymptotic to $g(\varphi)$, find an extra condition (if necessary) so that the uniform boundedness of (\ref{e1.2}) implies the existence of a compact uniform attractor.

\begin{definition}\label{Def2}
\cite{ArKo} \textup{
A square matrix $P$ with nonpositive off-diagonal entries is called
an M-matrix if one of the following equivalent conditions is met:
\begin{itemize}
\item[(a)] The leading principal minor determinants of $P$ are all positive.
\item[(b)] There is a vector $x>0$ (i.e. $x\in \textup{int}{\R}^N_+$) such that $Px >0$.
\item[(c)] There is a vector $y>0$ such that $P^Ty >0$.
\item[(d)] The matrix $P$ is nonsingular and the entries of $P^{-1}$ are all nonnegative.
\item[(e)] The real parts of the eigenvalues of $P$ are all positive, i.e. the matrix $-P$ is stable.
\end{itemize}}
\end{definition}

For any vector $u\in {\R}^N$, let $D(u) = \textup{diag}[u_1, \ldots, u_N]$.

\begin{theorem}\label{The2.1}
Assume that (\ref{e1.2}) meets the following requirements.
\begin{itemize}
\item[(i)] The $f_i$ are bounded on ${\R}_0\times S$ for any bounded set $S\subset C^+$.
\item[(ii)] For all $(t, \varphi)\in{\R}_0\times C^+$ and $i\in I_N$,
\begin{equation}\label{e2.1}
f_i(t, \varphi)\leq \beta_i + \sum^N_{j=1}a_{ij}\int^0_{-\tau}\varphi_j(\theta)d\xi_{ij}(\theta)-c_i\varphi_i(0),
\end{equation}
where $\beta_i>0$, $a_{ij}\geq 0$, $c_i>0$ and the $\xi_{ij}$ are nondecreasing with
\begin{equation}\label{e2.2}
\forall i, j\in I_N, \; \xi_{ij}(0)-\xi_{ij}(-\tau) = 1.
\end{equation}
\item[(iii)] The matrix $D(c)-A$ with $A=(a_{ij})$ is an M-matrix.
\end{itemize}
Then (\ref{e1.2}) has a compact uniform attractor.
\end{theorem}

\begin{remark}\label{Rem1} \textup{
The conditions of this theorem are irrelevant to either the size $\tau>0$ or the distribution of the delays. This eminent feature applies to all the results given in this section.}
\end{remark}
\begin{remark}\label{Rem2} \textup{
Conditions (iii) is crucial. When (iii) is not met, the conclusion may not be true; even the boundedness of solutions may no longer hold. For example, suppose $r_i>0$, $a_{ij}\geq 0$ for $i, j\in I_N \;(i\not=j)$ in (\ref{e1.1}). It is shown that \cite[Lemma 15.1.2]{HoSi} if $A$ has an eigenvalue $\lambda>0$ and a row vector $v\geq 0$ such that $vA = \lambda v$ then (\ref{e1.1}) has an unbounded solution in $\textup{int}{\R}^N_+$. Indeed \cite[Proof of Lemma 1]{Ho4}, (\ref{e1.1}) has an unbounded solution in $\textup{int}{\R}^N_+$ if $A$ is unstable, i.e. $-A$ is not an M-matrix.}
\end{remark}

\begin{theorem}\label{The2.2}
Assume that (\ref{e1.2}) meets the following requirements.
\begin{itemize}
\item[(i)] The $f_i$ are bounded on ${\R}_0\times S$ for any bounded set $S\subset C^+$.
\item[(ii)] For all $(t, \varphi)\in{\R}_0\times C^+$ and $i\in I_N$,
\begin{equation}\label{e2.3}
f_i(t, \varphi)\leq \beta_i - c_i\int^0_{-\tau}\varphi_i(\theta)d\xi_{ii}(\theta),
\end{equation}
where $\beta_i>0$, $c_i>0$ and the $\xi_{ii}$ are nondecreasing and satisfy (\ref{e2.2}).
\end{itemize}
Then (\ref{e1.2}) has a compact uniform attractor.
\end{theorem}
\begin{remark}\label{Rem3} \textup{
For $N=1$, Yang \cite{Ya} derived boundedness of solutions under some assumptions which are met if (ii) holds. However, from Definition \ref{Def1} we see that the property of having a compact uniform attractor is more that just boundedness.}
\end{remark}

The next result is the combination of Theorems \ref{The2.1} and \ref{The2.2} when the system can be arranged into triangular form of subsystems, of which each satisfies either Theorem \ref{The2.1} or \ref{The2.2}.

\begin{theorem}\label{The2.3}
Assume that $I_N$ has a partition $\{I^1, \ldots, I^m\}$ $(m >1)$ such that for $(t, \varphi)\in {\R}_0\times C^+$ with $\varphi = (\varphi^1, \ldots, \varphi^m)$,
\begin{eqnarray}
\forall i\in I^1, \; f_i(t, \varphi) &\leq& G^1_i(\varphi^1), \label{e2.4}\\
\forall i\in I^k (k>1), \; f_i(t, \varphi) &\leq& F^k_i(t, \varphi^1, \ldots, \varphi^{k-1})+G^k(\varphi^k), \label{e2.5}
\end{eqnarray}
where the $F^k$ are bounded on ${\R}_0\times S$ for any bounded set $S\subset C^+$ and each $G^k_i$ has either the form
\begin{equation}\label{e2.6}
\forall i\in I^k,\; G^k_i(\varphi^k)= \beta_i + \sum_{j\in I^k}a_{ij}\int^0_{-\tau}\varphi^k_j(\theta)d\xi_{ij}(\theta)-c_i\varphi^k_i(0),
\end{equation}
with $D(c)^k-A^k$ an M-matrix, or the form
\begin{equation}\label{e2.7}
\forall i\in I^k,\; G^k_i(\varphi^k)= \beta_i -c_i\int^0_{-\tau}\varphi^k_i(\theta)d\xi_{ii}(\theta),
\end{equation}
where the $\beta_i$, $a_{ij}$, $c_i$ and $\xi_{ij}$ are the same as in (\ref{e2.1}) and $D(c)^k-A^k$ is the corresponding $|I^k|\times |I^k|$ matrix. Then (\ref{e1.2}) has a compact uniform attractor.
\end{theorem}

The main purpose of establishing these theorems for a system to have a compact uniform attractor is to apply this property to the study of permanence.

\begin{definition}\label{Def3} \textup{
For any nonempty set $J\subset I_N$, (\ref{e1.2}) is said to be partially permanent with respect to $J$ if there exist $\delta>0$ and $M>\delta$ such that, for all $(t_0, \varphi)\in {\R}_0\times C^+$ with $\varphi_i(0)>0$ for all $i\in J$, the solution of (\ref{e1.2}) satisfies
\[
\forall i\in J, \forall\; \textup{large}\; t, \; \delta < x_i(t, t_0, \varphi) < M.
\]}
\end{definition}
From this definition we see that (\ref{e1.2}) is permanent if it is partially permanent with respect to $J = I_N$.

For any $i\in I_N$ and $\Omega\subset C^+$, let $\Omega_i = \{\varphi\in\Omega: \varphi_i =0\}$.

\begin{theorem}\label{The2.4}
Assume that (\ref{e1.2}) satisfies the following conditions.
\begin{itemize}
\item[(i)] (\ref{e1.2}) has a compact uniform attractor $\Omega\subset C^+$.
\item[(ii)] $f$ is bounded on ${\R}_0\times \Omega$ and uniformly Lipschitzian on $\Omega$, i.e. there is a $K>0$ such that
    \[
    \forall t\in {\R}_0, \forall \varphi,\psi\in\Omega, \; |f(t, \varphi)-f(t, \psi)| \leq K\|\varphi-\psi\|.
    \]
\item[(iii)] For a nonempty set $J\subset I_N$ and $\varphi\in \cup_{i\in J}\Omega_i$, $f(t, \varphi)$ is $T_0$-periodic.
\item[(iv)] There are $q_i>0$, $i\in J$, for each $(t_0, \varphi)\in {\R}_0\times (\cup_{i\in J}\Omega_i)$, there is a $T(t_0, \varphi)>0$ such that
    \begin{equation}\label{e2.8}
    \int^{T(t_0, \varphi)}_0\sum_{i\in J}q_if_i(t_0+s, x_{t_0+s}(t_0, \varphi))ds >0.
    \end{equation}
\end{itemize}
Then (\ref{e1.2}) is partially permanent with respect to $J$. If $J= I_N$ then (\ref{e1.2}) is permanent.
\end{theorem}

\begin{remark}\label{Rem4} \textup{
This theorem when $J=I_N$ can be viewed as an extension of \cite[Theorem 12.2.1]{HoSi} from a system without delays on a closed set $S_n\subset {\R}^N_+$ to (\ref{e1.2}) on ${\R}_0\times C^+$. Unfortunately, it is not easily applicable to any concrete system as condition (iv) is hardly checkable. However, we can develop an easily checkable condition for a class of systems based on this.}
\end{remark}

A particular case of (\ref{e1.2}) is that
\begin{eqnarray}
f_i(t, \varphi) &=& r_i(t) + L_i(\varphi) - F_i(t, \varphi), \label{e2.9}\\
L_i(\varphi) &=& \sum^N_{j=1}a_{ij}\int^0_{-\tau}\varphi_j(\theta)d\xi_{ij}(\theta) - \sum^N_{j=1}b_{ij}\int^0_{-\tau}\varphi_j(\theta)d\eta_{ij}(\theta), \label{e2.10}\\
\forall i, j\in I_N, && a_{ij}\geq 0, b_{ij}\geq 0, b_{ii}>0, \label{e2.11}
\end{eqnarray}
where the $r_i$ are continuous $T_0$-periodic with
\begin{equation}\label{e2.12}
\forall i\in I_N, \forall t_0\in {\R}_0, \; \frac{1}{T_0}\int^{T_0}_0r_i(t_0+s)ds = \bar{r}_i > 0,
\end{equation}
and the $\xi_{ij}$ and $\eta_{ij}$ are nondecreasing with
\begin{equation}\label{e2.13}
\forall i,j\in I_N, \; \xi_{ij}(0)-\xi_{ij}(-\tau) = 1, \eta_{ij}(0)-\eta_{ij}(-\tau) = 1.
\end{equation}
We assume that the $F_i$ are nonnegative, bounded on ${\R}_0\times S$ and uniformly Lipschitzian in $\varphi\in S$ for any bounded set $S\subset C^+$. We assume also that (\ref{e1.2}) with (\ref{e2.9})--(\ref{e2.13}) meets the requirement of one of Theorems \ref{The2.1}--\ref{The2.3} or Proposition \ref{Pro1}. Then (\ref{e1.2}) with (\ref{e2.9})--(\ref{e2.13}) has a compact uniform attractor $\Omega\subset C^+$ and $f$ is bounded on ${\R}_0\times\Omega$ and uniformly Lipschitzian on $\Omega$. From now on (\ref{e1.2}) with (\ref{e2.9})--(\ref{e2.13}) is always assumed to have these properties.

In addition to (\ref{e1.2}) with (\ref{e2.9})--(\ref{e2.13}), consider also the autonomous Lotka-Volterra system
\begin{equation}\label{e2.14}
x'_i = x_i(\bar{r}_i +(A-B)_ix), \quad i\in I_N,
\end{equation}
where $A=(a_{ij})$ and $B=(b_{ij})$. Denote the $i$th coordinate plane by $\pi_i = \{x\in{\R}^N_+: x_i=0\}$.

\begin{theorem}\label{The2.5}
For (\ref{e1.2}) with (\ref{e2.9})--(\ref{e2.13}) and a nonempty set $J\subset I_N$, assume also that
\begin{equation}\label{e2.15}
\forall j\in I_N, \forall t\in {\R}_0, \forall \varphi \in \cup_{i\in J}\Omega_i, \; F_j(t, \varphi)\equiv 0,
\end{equation}
\begin{equation}\label{e2.16}
\exists q_i>0\; \textup{for}\; i\in J\; \textup{such that}\; \sum_{i\in J}q_i(\bar{r}_i+(A-B)_i\hat{x})>0
\end{equation}
for all fixed point $\hat{x}$ of (\ref{e2.14}) in $\cup_{i\in J}\pi_i$. Then (\ref{e1.2}) with (\ref{e2.9})--(\ref{e2.13}) is partially permanent with respect to $J$. If also $J= I_N$ then (\ref{e1.2}) with (\ref{e2.9})--(\ref{e2.13}) is permanent.
\end{theorem}

\begin{remark}\label{Rem5} \textup{
This theorem when $J=I_N$ is the extension of Jansen's result \cite{Ja} from autonomous replicator and Lotka Volterra systems to (\ref{e1.2}) with (\ref{e2.9})--(\ref{e2.13}) (see also \cite[Theorem 13.6.1 and Exercise 13.6.3]{HoSi}).}
\end{remark}

Note that permanence of (\ref{e1.2}) with respect to $\{i\}$ for every $i\in J$ implies permanence with respect to $J$. Then applying Theorems \ref{The2.4} and \ref{The2.5} to each $i\in J$ we obtain the following corollaries.

\begin{corollary}\label{Cor1}
Assume that (\ref{e1.2}) satisfies the conditions (i)--(iii) of Theorem \ref{The2.4}. Moreover, for each $i\in J$ and every $(t_0, \varphi)\in{\R}_0\times \Omega_i$, there is a $T(t_0, \varphi)>0$ such that
\begin{equation}\label{e2.17}
\int^{T(t_0, \varphi)}_0f_i(t_0+s, x_{t_0+s}(t_0, \varphi))ds >0.
\end{equation}
Then (\ref{e1.2}) is partially permanent with respect to $J$. If also $J= I_N$ then (\ref{e1.2}) is permanent.
\end{corollary}

\begin{corollary}\label{Cor2}
Assume that (\ref{e1.2}) with (\ref{e2.9})--(\ref{e2.13}) satisfies (\ref{e2.15}). Moreover, for each $i\in J$ and every fixed point $\hat{x}$ of (\ref{e2.14}) in $\pi_i$, we have $\bar{r}_i+ (A-B)_i\hat{x}>0$. Then (\ref{e1.2}) with (\ref{e2.9})--(\ref{e2.13}) is partially permanent with respect to $J$. If also $J= I_N$ then (\ref{e1.2}) with (\ref{e2.9})--(\ref{e2.13}) is permanent.
\end{corollary}

When $J=I_N$, if we apply Corollary \ref{Cor2} to every subsystem of (\ref{e1.2}) with (\ref{e2.9})--(\ref{e2.13}) and (\ref{e2.15}), we obtain the following.

\begin{corollary}\label{Cor3}
Assume that (\ref{e1.2}) with (\ref{e2.9})--(\ref{e2.13}) and (\ref{e2.15}) satisfies
\begin{equation}\label{e2.18}
\forall i\in I_N, \forall\hat{x}\in \pi_i\ \;\textup{(fixed points of (\ref{e2.14}))},\; \bar{r}_i+ (A-B)_i\hat{x}>0.
\end{equation}
Then (\ref{e1.2}) and all of its subsystems are permanent.
\end{corollary}
\begin{remark}\label{Rem6} \textup{
Mierczy\'{n}ski and Schreiber \cite[Corollary 3.1]{MiSc} proved that if (\ref{e2.14}) is dissipative then (\ref{e2.18}) is a necessary and sufficient condition for (\ref{e2.14}) and all of its subsystems to be robustly permanent, i.e. all small perturbations of the systems (in some sense) are permanent. Ahmad and Lazer \cite{AhLa4} extended this result to nonautonomous competitive Lotka-Volterra system
\begin{equation}\label{e2.19}
x'_i(t) = x_i(t)[r_i(t) - A_ix(t)], \quad i\in I_N,
\end{equation}
where the $r_i$ are bounded continuous satisfying
\begin{equation}\label{e2.20}
\forall i\in I_N, \lim_{T\to\infty}m(r_i, t_0, t_0+T) = \bar{r}_i>0 \; \textup{uniformly for}\; t_0\in {\R}_0.
\end{equation}
Hou \cite{Ho3} and \cite{Ho4} further extended it to delayed nonautonomous Lotka-Volterra systems of the form
\begin{equation}\label{e2.21}
x'_i(t) = x_i(t)[r_i(t) + L_i(x_t)], \quad, i\in I_N,
\end{equation}
where the $r_i$ are the same as in (\ref{e2.19}) and the $L_i(\varphi)$ are given by (\ref{e2.10}), (\ref{e2.11}) and (\ref{e2.13}). Here Corollary \ref{Cor3} can be viewed as a further partial extension of \cite[Corollary 3.1]{MiSc} to (\ref{e1.2}) with (\ref{e2.9})--(\ref{e2.13}) and (\ref{e2.15}).}
\end{remark}

\noindent \textbf{Open Problem} Can Theorem \ref{The2.5} and Corollaries \ref{Cor2} and \ref{Cor3} be extended to (\ref{e1.2}) with (\ref{e2.9})--(\ref{e2.11}), (\ref{e2.13}) and (\ref{e2.15}) but the $r_i$ satisfy (\ref{e2.20}) instead of being $T_0$-periodic?

\section{Proof of the existence of a compact uniform attractor}

\textbf{Proof of Theorem \ref{The2.1}}
Since $D(c)-A$ is an M-matrix, there is a vector $d\in\textrm{int}{\R}^N_+$ such that $(D(c)-A)d >0$, i.e.
\begin{equation}\label{e3.1}
\forall i\in I_N, \; c_id_i-A_id >0.
\end{equation}
Let
\begin{equation}
M_0 = \max\biggl\{\frac{\beta_i}{c_id_i-A_id}: i\in I_N\biggr\}, \label{e3.2}
\end{equation}
\begin{equation}
\forall \varphi \in  C^+, \; \|\varphi\|_{d^{-1}} = \sup_{\theta\in[-\tau, 0]}\biggl\{\max_{i\in I_N} \varphi_i(\theta)d^{-1}_i\biggr\}. \label{e3.3}
\end{equation}
We divide the rest of the proof into the following five steps.

\textbf{Step 1} For any $\varphi\in C^+$ with $\varphi(0)\not= 0$ and $t_0\in{\R}_0$, the solution $x_t(t_0, \varphi)$ exists on $[t_0, \infty)$ and satisfies
\begin{equation}\label{e3.4}
\|x_t(t_0, \varphi)\|_{d^{-1}} \leq \max\{\|\varphi\|_{d^{-1}}, M_0\}.
\end{equation}
For if there is a $t_1 >t_0$ such that $\|x_{t_1}\|_{d^{-1}}>
\max\{\|\varphi\|_{d^{-1}}, M_0\}$, then there are $t_2 \in (t_0,
t_1]$ and $i\in I_N$ satisfying
\[
x_i(t_2)d_i^{-1} = |x(t_2)|_{d^{-1}}=\max_{t_0-\tau\leq t\leq
t_1}|x(t)|_{d^{-1}} = \|x_{t_2}\|_{d^{-1}} \geq
\|x_{t_1}\|_{d^{-1}} > M_0.
\]
This implies $x'_i(t_2) \geq 0$. On the other hand, however, from
(\ref{e1.2}), (\ref{e2.1}) and (\ref{e3.1})--(\ref{e3.3}) we have
\begin{eqnarray*}
x'_i(t_2) &\leq& x_i(t_2)[\beta_i + A_id\|x_{t_2}\|_{d^{-1}} -
 c_id_ix_i(t_2)d_i^{-1}]\\
 &=& x_i(t_2)[\beta_i - (c_id_i-A_id)\|x_{t_2}\|_{d^{-1}}] \\
 &\leq & x_i(t_2)(c_id_i-A_id)(M_0-\|x_{t_2}\|_{d^{-1}}) < 0.
\end{eqnarray*}
This contradiction shows the truth of (\ref{e3.4}) on the existing
interval of $x$ and (\ref{e3.4}) ensures the extension of the
solution to $[t_0, \infty)$.

\textbf{Step 2} We claim that if $\|\varphi\|_{d^{-1}}> M_0$ then $\|x_t\|_{d^{-1}}$ is nonincreasing as long as $\|x_t\|_{d^{-1}}\geq M_0$. For if there are $t_2 >t_1 \geq t_0$ such that
\[
\|x_{t_2}\|_{d^{-1}} > \|x_{t_1}\|_{d^{-1}} > M_0,
\]
then there are $t_3\in (t_1, t_2]$ and $i\in I_N$ satisfying
\[
x_i(t_3)d^{-1}_i = |x(t_3)|_{d^{-1}} = \max_{t_1-\tau \leq t\leq t_2}|x(t)|_{d^{-1}} = \|x_{t_3}\|_{d^{-1}} \geq \|x_{t_2}\|_{d^{-1}}.
\]
Thus, $x'_i(t_3)\geq 0$. However, by the same technique as that used in the proof of (\ref{e3.4}), we derive $x'_i(t_3)<0$, a contradiction to $x'_i(t_3)\geq 0$. This shows the truth of our claim.

\textbf{Step 3} We show that for every $(t_0, \varphi)\in {\R}_0\times C^+$,
\begin{equation}\label{e3.5}
\limsup_{t\to\infty}\|x_t(t_0, \varphi)\|_{d^{-1}} \leq M_0.
\end{equation}
If this is not true, then some solution satisfies $\limsup_{t\to\infty} \|x_t\|_{d^{-1}} > M_0$. From step 2 we know that $\|x_t\|_{d^{-1}}$ is nonincreasing so $\lim_{t\to\infty} \|x_t\|_{d^{-1}} =\bar{M}_0 > M_0$. Let
\[
u_i = \limsup_{t\to\infty} x_i(t)d_i^{-1}, \quad u_j = \max\{u_i: i\in I_N\}
\]
for all $i$ and some $j$ in $I_N$. We look for an increasing sequence $\{t_k\}$ with $t_k\to\infty$ as $k\to\infty$ such that
\begin{equation}
\lim_{k\to\infty} x'_j(t_k) = 0 \;\;\textup{and}\;\; \lim_{k\to\infty}x_j(t_k)d_j^{-1} = u_j. \label{e3.6}
\end{equation}
If $x_j(t)$ is not monotone for large $t$ then we take a sequence $\{t_k\}$ so that each $x_j(t_k)$ is a local maximum of $x_j(t)$. This sequence certainly fulfils (\ref{e3.6}). If $x_j(t)$ is nondecreasing (nonincreasing) then $x'_j(t)\geq 0 \; (\leq 0)$ and, by the boundedness of $x$, $\liminf_{t\to\infty}x'_j(t) = 0$ ($\limsup_{t\to\infty}x'_j(t) = 0$). Then we can choose a sequence $\{t_k\}$ satisfying (\ref{e3.6}).

We check that $u_j = \bar{M}_0$. From (\ref{e3.6}) we have
\[
u_j = \lim_{k\to\infty}x_j(t_k)d_j^{-1} \leq \lim_{k\to\infty} \|x_{t_k}\|_{d^{-1}} = \lim_{t\to\infty}\|x_t\|_{d^{-1}} = \bar{M}_0.
\]
On the other hand, for any $\varepsilon >0$, by the definition of $u_j$ there exists $T>t_0$ such that
\[
\forall i\in I_N, \; \forall t\geq T, \;\; x_i(t)d_i^{-1} < u_j +\varepsilon.
\]
From this follows $\|x_t\|_{d^{-1}} < u_j + \varepsilon$ for $t\geq T+\tau$. Since $\|x_t\|_{d^{-1}}$ is nonincreasing, we have $\bar{M}_0 < u_j+\varepsilon$ so $\bar{M}_0 \leq u_j$ as $\varepsilon\to 0+$. Therefore, $u_j = \bar{M}_0$. Then it follows from this and (\ref{e3.6}) that
\[
\lim_{k\to\infty}\varepsilon_k = \lim_{k\to\infty} \delta_k = 0,
\]
where
\[
\forall k\geq 1, \;\;\varepsilon_k = x_j(t_k)d_j^{-1}-\bar{M}_0, \;\; \delta_k = \|x_{t_k}\|_{d^{-1}} - \bar{M}_0.
\]
Thus, from (\ref{e1.2}), (\ref{e2.1}),  (\ref{e3.1}) and (\ref{e3.2}),
\begin{eqnarray*}
x'_j(t_k) &\leq & x_j(t_k)[\beta_j + A_jd\|x_{t_k}\|_{d^{-1}} - c_jd_jx_j(t_k)d_j^{-1}] \\
 &=& x_j(t_k)[\beta_j +A_jd(\bar{M}_0+\delta_k) -c_jd_j(\bar{M}_0+\varepsilon_k)] \\
 &\leq & x_j(t_k)(c_jd_j-A_jd)\biggl[M_0-\bar{M}_0 +\frac{A_jd\delta_k-c_jd_j\varepsilon_k}{c_jd_j-A_jd}\biggr] \\
 &\to & \bar{M}_0d_j(c_jd_j-A_jd)(M_0-\bar{M}_0) \quad (k\to\infty) \\
 & <& 0.
\end{eqnarray*}
This contradicts (\ref{e3.6}) and hence shows (\ref{e3.5}).

\textbf{Step 4} We show that for each $M_1 >M_0$ and every $\varphi\in C^+$ with $\varphi(0)\not= 0$, there is a $T = T(\varphi, M_1)>0$ such that
\begin{equation}
\forall t_0\in {\R}_0, \; \forall t\geq t_0+T, \;\; \|x_t(t_0, \varphi)\|_{d^{-1}} < M_1. \label{e3.7}
\end{equation}
Suppose this is not true. Then, for some $\varphi\in C^+$ with $\varphi(0)\not= 0$ and $M_1>M_0$, by (\ref{e3.5}) there are $\{t_k\}\subset {\R}_0$ and $\{T_k\}\subset [3\tau, \infty)$, $T_k\uparrow\infty$ as $k\to\infty$, such that
\begin{equation}
\forall k\geq 1,  \|x_{t_k+T_k}(t_k, \varphi)\|_{d^{-1}} = M_1, \; \forall t> t_k+T_k,  \|x_t(t_k, \varphi)\|_{d^{-1}} <M_1. \label{e3.8}
\end{equation}
From step 2, $\|x_t(t_k, \varphi)\|_{d^{-1}}$ is nonincreasing for $t\in [t_k, t_k+T_k]$. Since the function $\|x_t(t_k, \varphi)\|_{d^{-1}}$ is continuous, with
\[
g_k(t) = \|x_{t-2\tau}(t_k, \varphi)\|_{d^{-1}} - \|x_t(t_k, \varphi)\|_{d^{-1}},
\]
$g_k$ is continuous and nonnegative for $t\in [t_k+2\tau, t_k+T_k]$. So there is an $s_k\in [t_k+2\tau, t_k+T_k]$ such that
\[
g_k(s_k) = \min\{g_k(t): t\in [t_k+2\tau, t_k+T_k]\}.
\]

(4a) We first show that
\begin{equation}
\lim_{k\to\infty}g_k(s_k) = 0. \label{e3.9}
\end{equation}
Suppose (\ref{e3.9}) is not true. Then $\limsup_{k\to\infty}g_k(s_k) >0$. By choosing a subsequence of $\{s_k\}$ if necessary, we may assume the existence of $p_0>0$ and an integer $K$ such that, for all $k\geq K$, $g_k(s_k)\geq p_0$ so that
\begin{equation}
\forall t\in [t_k+2\tau, t_k+T_k], \; \|x_{t-2\tau}(t_k, \varphi)\|_{d^{-1}} \geq \|x_t(t_k, \varphi)\|_{d^{-1}} +p_0. \label{e3.10}
\end{equation}
Let $m_k$ be the greatest integer part of $T_k/(2\tau)$. Then $\{m_k\}$ is unbounded due to the unboundedness of $\{T_k\}$. However, by step 2 and repeatedly using (\ref{e3.10}), we have
\begin{eqnarray*}
\|\varphi\|_{d^{-1}} &\geq & \|x_{t_k+T_k-2\tau m_k}(t_k, \varphi)\|_{d^{-1}} \\
 &\geq & \|x_{t_k+T_k-2\tau(m_k-1)}(t_k, \varphi)\|_{d^{-1}} + p_0 \\
 &\geq & m_kp_0 +\|x_{t_k+T_k}(t_k, \varphi)\|_{d^{-1}} \\
 &=& m_kp_0 +M_1.
\end{eqnarray*}
This contradiction to the unboundedness of $\{m_k\}$ shows the truth of (\ref{e3.9}).

(4b) We next show that for each $k\geq 1$, there is an $\ell_k\in [s_k-\tau, s_k]$ satisfying
\begin{equation}
\biggl(\frac{d}{dt}\|x_t(t_k, \varphi)\|_{d^{-1}}\biggr)_{t=\ell_k} \geq -\frac{2}{\tau}g_k(s_k). \label{e3.11}
\end{equation}
If $g_k(s_k) =0$ then
\[
\|x_t(t_k, \varphi)\|_{d^{-1}}\equiv \|x_{s_k}(t_k, \varphi)\|_{d^{-1}}\; \textup{for}\; t\in [s_k-2\tau, s_k]
\]
so $\frac{d}{dt}\|x_t(t_k, \varphi)\|_{d^{-1}} = 0$ for $t\in (s_k-2\tau, s_k)$. In this case, $\ell_k=s_k-\tau$ meets the requirement of (\ref{e3.11}). Suppose $g_k(s_k)>0$ and (\ref{e3.11}) does not hold for any $\ell_k\in [s_k-\tau, s_k]$. Then $\frac{d}{dt}\|x_t(t_k, \varphi)\|_{d^{-1}} < -\frac{2}{\tau}g_k(s_k)$ for almost every $t\in [s_k-\tau, s_k]$ so
\begin{eqnarray*}
g_k(s_k) &=& \|x_{s_k-2\tau}(t_k, \varphi)\|_{d^{-1}} - \|x_{s_k}(t_k, \varphi)\|_{d^{-1}} \\
 &\geq & \|x_{s_k-\tau}(t_k, \varphi)\|_{d^{-1}} -\|x_{s_k}(t_k, \varphi)\|_{d^{-1}} \\
 &=& -\int^{s_k}_{s_k-\tau}\biggl(\frac{d}{dt}\|x_t(t_k, \varphi)\|_{d^{-1}}\biggr)dt \\
 &\geq & 2 g_k(s_k).
\end{eqnarray*}
This contradiction to $0<g_k(s_k) <2g_k(s_k)$ shows the existence of $\ell_k$ satisfying (\ref{e3.11}).

(4c) We further show that for each $k\geq 1$, there are $w_k\in[\ell_k-\tau, \ell_k]$ and $i_k\in I_N$ such that
\begin{equation}
x'_{i_k}(w_k, t_k, \varphi) \geq -\frac{2d_{i_k}}{\tau}g_k(s_k). \label{e3.12}
\end{equation}
Indeed, for each $k\geq 1$, there are $w_k\in [\ell_k-\tau, \ell_k]$ and $i_k\in I_N$ such that
\[
\|x_{\ell_k}(t_k, \varphi)\|_{d^{-1}} = |x(w_k, t_k, \varphi)|_{d^{-1}} = x_{i_k}(w_k, t_k, \varphi)d^{-1}_{i_k}.
\]
For small $\delta >0$, $w_k-\delta\in [\ell_k-\tau-\delta, \ell_k-\delta]$ so
\begin{eqnarray*}
\qquad & & -\frac{1}{\delta}[x_{i_k}(w_k-\delta, t_k, \varphi)-x_{i_k}(w_k, t_k, \varphi)]d^{-1}_{i_k} \\
 &\geq & -\frac{1}{\delta}[|x(w_k-\delta, t_k, \varphi)|_{d^{-1}} - \|x_{\ell_k}(t_k, \varphi)\|_{d^{-1}}] \\
 &\geq & -\frac{1}{\delta}[\|x_{\ell_k-\delta}(t_k, \varphi)\|_{d^{-1}} - \|x_{\ell_k}(t_k, \varphi)\|_{d^{-1}}].
\end{eqnarray*}
As $\delta\to 0+$, the above inequalities lead to
\[
x'_{i_k}(w_k, t_, \varphi)d^{-1}_{i_k} \geq \biggl(\frac{d}{dt}\|x_t(t_k, \varphi)\|_{d^{-1}}\biggr)_{t=\ell_k}.
\]
Then (\ref{e3.12}) follows from this and (\ref{e3.11}).

Now armed with (4a)--(4c) we are able to construct a contradiction. By choosing a subsequence of $\{k\}$ if necessary, without loss of generality, we may assume that $i_k=i_0\in I_N$ for all $k\geq 1$ so that (\ref{e3.12}) becomes
\[
x'_{i_0}(w_k, t_k, \varphi) \geq -\frac{2d_{i_0}}{\tau}g_k(s_k).
\]
It then follows from this and (\ref{e3.9}) that
\begin{equation}
\liminf_{k\to\infty}x'_{i_0}(w_k, t_k, \varphi) \geq 0. \label{e3.13}
\end{equation}
Nevertheless, from (\ref{e1.2}), (\ref{e2.1}), (\ref{e3.1}), (\ref{e3.2}), step 2 and the equalities below (\ref{e3.12}),
\begin{eqnarray*}
\quad & & x'_{i_0}(w_k, t_k, \varphi)/x_{i_0}(w_k, t_k, \varphi) \\
&\leq & \beta_{i_0} +A_{i_0}d\|x_{w_k}(t_k, \varphi)\|_{d^{-1}} -c_{i_0}d_{i_0}x_{i_0}(w_k, t_k, \varphi)d^{-1}_{i_0} \\
&\leq & \beta_{i_0} +A_{i_0}d\|x_{s_k-2\tau}(t_k, \varphi)\|_{d^{-1}} -c_{i_0}d_{i_0}\|x_{\ell_k}(t_k, \varphi)\|_{d^{-1}} \\
&\leq & \beta_{i_0} +A_{i_0}d(\|x_{s_k}(t_k, \varphi)\|_{d^{-1}}+g_k(s_k)) -c_{i_0}d_{i_0}\|x_{s_k}(t_k, \varphi)\|_{d^{-1}} \\
&\leq & \beta_{i_0} -(c_{i_0}d_{i_0}-A_{i_0}d)M_1 +A_{i_0}dg_k(s_k) \\
&\leq & (c_{i_0}d_{i_0}-A_{i_0}d)(M_0-M_1)+A_{i_0}d g_k(s_k).
\end{eqnarray*}
From (\ref{e3.8}) and step 2 we know that
\[
x_{i_0}(w_k, t_k, \varphi)d^{-1}_{i_0} = \|x_{\ell_k}(t_k, \varphi)\|_{d^{-1}} \geq \|x_{s_k}(t_k, \varphi)\|_{d^{-1}} \geq M_1.
\]
As $M_1 >M_0$, from (\ref{e3.9}) and the above inequalities we obtain
\[
\limsup_{k\to\infty}x'_{i_0}(w_k, t_k, \varphi) \leq d_{i_0}M_1(c_{i_0}d_{i_0}-A_{i_0}d)(M_0-M_1)< 0.
\]
This contradiction to (\ref{e3.13}) shows the existence of $T = T(\varphi, M_1)>0$ satisfying (\ref{e3.7}).

\textbf{Step 5} Construction of a compact uniform attractor $\Omega\subset C^+$. For this purpose, fix an $M_1 > M_0$ and let
\[
S_0 = \{\varphi\in C^+: \|\varphi\|_{d^{-1}}\leq M_1\}.
\]
Then, by the assumption on $f$, there is a $\rho>0$ such that
\[
\forall i\in I_N, \forall (t, \varphi)\in {\R}_0\times S_0, \; |f_i(t, \varphi)|\leq \rho.
\]
Define
\begin{equation}\label{e3.14}
\Omega = \{\varphi\in S_0: \forall i\in I_N, \forall \theta_1,\theta_2\in [-\tau, 0]\; \textup{with}\; \theta_1 \leq \theta_2, (\ref{e3.15})\; \textup{holds}\},
\end{equation}
\begin{equation}\label{e3.15}
\varphi_i(\theta_1) e^{-\rho(\theta_2-\theta_1)} \leq \varphi_i(\theta_2 \leq \varphi_i(\theta_1) e^{\rho(\theta_2-\theta_1)}.
\end{equation}
We check that this $\Omega$ is a compact uniform attractor.

(i) For any convergent sequence $\{\varphi^n\}\subset \Omega$ with a limit $\varphi^0\in C^+$, since all the inequalities for each $\varphi^n$ in the definition of $S_0$ and $\Omega$ are retained for $\varphi^0$, we have $\varphi^0\in \Omega$ so $\Omega$ is closed. By (\ref{e3.15}) we have
\[
\varphi_i(\theta_1)[e^{-\rho(\theta_2-\theta_1)}-1] \leq \varphi_i(\theta_2)-\varphi_1(\theta_1) \leq \varphi_i(\theta_1)[e^{\rho(\theta_2-\theta_1)}-1]
\]
so
\[
|\varphi_i(\theta_2)-\varphi_i(\theta_1)| \leq |d|M_1\rho e^{\tau\rho}|\theta_2-\theta_1|.
\]
This shows that the functions in $\Omega$ over $[-\tau, 0]$ are equicontinuous. Since $\Omega$ is bounded, by Arzela-Ascoli theorem, $\Omega$ is relatively compact. This, together with the closedness, shows that $\Omega$ is compact.

(ii) For each $\varphi\in \Omega$, from step 2 we know that $x_t(t_0, \varphi)\in S_0$ for all $t_0\in{\R}_0$ and $t\geq t_0$. For any $t_2 \geq t_1 \geq t_0$, integration of (\ref{e1.2}) gives
\[
x_i(t_2, t_0, \varphi) = x_i(t_1, t_0, \varphi)\exp\biggl(\int^{t_2}_{t_1}f_i(s, x_s(t_0, \varphi))ds\biggr)
\]
so
\[
x_i(t_1, t_0, \varphi)e^{-\rho(t_2-t_1)} \leq x_i(t_2, t_0, \varphi) \leq x_i(t_1, t_0, \varphi)e^{\rho(t_2-t_1)}.
\]
This shows that $x_t(t_0, \varphi)\in \Omega \not=\emptyset$ for $t\geq t_0 +\tau$. For $t\in (t_0, t_0+\tau)$ and $-\tau \leq \theta_1 < \theta_2\leq 0$, if $t+\theta_1 \geq t_0$ then
\[
x_i(t+\theta_2, t_0, \varphi) = x_i(t+\theta_1, t_0, \varphi)\exp\biggl(\int^{\theta_2}_{\theta_1}f_i(t+s, x_{t+s}(t_0, \varphi))ds\biggr)
\]
so
\begin{equation}\label{e3.16}
x_i(t+\theta_1, t_0, \varphi)e^{-\rho(\theta_2-\theta_1)} \leq x_i(t+\theta_2, t_0, \varphi) \leq x_i(t+\theta_1, t_0, \varphi)e^{\rho(\theta_2-\theta_1)}.
\end{equation}
If $t+\theta_2\leq t_0$ then $x_t(\theta_j, t_0, \varphi) = \varphi(t-t_0+\theta_j)$ for $j = 1, 2$ so (\ref{e3.16}) follows from (\ref{e3.15}). If $t+\theta_1 < t_0 < t+\theta_2$, then
\begin{eqnarray*}
x_i(t+\theta_2, t_0, \varphi) &\leq & \varphi_i(0)e^{\rho(t+\theta_2-t_0)} \\
    &\leq & \varphi_i(t-t_0+\theta_1)e^{\rho(t_0-t-\theta_1)}e^{\rho(t+\theta_2-t_0)}\\
    &= & x_i(t+\theta_1, t_0, \varphi)e^{\rho(\theta_2-\theta_1)}, \\
x_i(t+\theta_2, t_0, \varphi) &\geq & \varphi_i(0)e^{-\rho(t+\theta_2-t_0)} \\
    &\geq & \varphi_i(t-t_0+\theta_1)e^{-\rho(t_0-t-\theta_1)}e^{-\rho(t+\theta_2-t_0)}\\
    &= & x_i(t+\theta_1, t_0, \varphi)e^{-\rho(\theta_2-\theta_1)}.
\end{eqnarray*}
Thus, (\ref{e3.16}) holds and $x_t(t_0, \varphi)\in \Omega$. Therefore, $\varphi\in\Omega$ implies $x_t(t_0, \varphi)\in \Omega$ for all $t_0\in{\R}_0$ and all $t\geq t_0$.

(iii) For each $\varphi\in C^+$, from step 4 we know the existence of $T = T(\varphi)>0$ such that $x_t(t_0, \varphi)\in S_0$ holds for all $t_0\in{\R}_0$ and $t\geq t_0+T$. Then from (ii) we obtain $x_t(t_0, \varphi)\in \Omega$ for $t\geq t_0+T+\tau$.

(iv) For each $\varphi\in\Omega$, (\ref{e3.15}) implies that $\varphi_i(\theta_0)=0$ for some $i\in I_N$ and some $\theta_0\in [-\tau, 0]$ if and only if $\varphi_i(\theta)\equiv 0$ on $[-\tau, 0]$.

Therefore, $\Omega$ defined by (\ref{e3.14}) is a compact uniform attractor of (\ref{e1.2}).
\rule{2mm}{3mm}

\textbf{Proof of Theorem \ref{The2.2}}
The techniques of the proof are the same as those of \cite[Lemmas 1, 3, 4]{Ho3} and part of \cite[Lemma 3.1]{TaKu}. But for completeness we give a full outline following the proof of Theorem \ref{The2.1}. For  $(t_0, \varphi)\in {\R}_0\times C^+$ with $\varphi_i(0)>0$, we have
\[
x'_i(t, t_0, \varphi) \leq \beta_ix_i(t, t_0, \varphi)
\]
so for $t\geq t_0+\tau$ and $\theta\in [-\tau, 0]$, $(x_i)_t(\theta)\geq x_i(t, t_0, \varphi)e^{\beta_i\theta}$. Hence, for $t\geq t_0+\tau$ in its existing interval,
\begin{equation}\label{e3.17}
x'_i(t, t_0, \varphi) \leq x_i(t, t_0, \varphi)\biggl[\beta_i-c_i\int^0_{-\tau}e^{\beta_i\theta}d\xi_{ii}(\theta)x_i(t, t_0, \varphi)\biggr].
\end{equation}
Let
\begin{equation}\label{e3.18}
\rho_i = \frac{\beta_i}{c_i\int^0_{-\tau}e^{\beta_i\theta}d\xi_{ii}(\theta)}, \quad \rho^0 = \max\{\rho_i: i\in I_N\}.
\end{equation}
Then $x_i(t, t_0, \varphi)$ is decreasing as long as $x_i(t, t_0, \varphi)>\rho_i$. This shows the existence and boundedness of $x(t, t_0, \varphi)$ on $[t_0, \infty)$.

Now multiplying (\ref{e3.17}) by $-x^{-2}_i(t, t_0, \varphi)e^{\beta_i(t-t_0)}$ and integrating, we obtain
\begin{equation}\label{e3.19}
x_i(t_0+t, t_0, \varphi) \leq \frac{e^{\beta_it}}{\varphi_i(0)^{-1}+\rho_i^{-1}(e^{\beta_it}-1)}.
\end{equation}
From this follows
\[
\limsup_{t\to\infty} x_i(t, t_0, \varphi) \leq \rho_i\leq \rho^0.
\]
Moreover, for any fixed $\rho >\rho^0$, (\ref{e3.19}) shows that for all $\varphi\in C^+$,  there exists a $T=T(\varphi)>0$ such that
\begin{equation}\label{e3.20}
\forall t_0\in{\R}_0, \forall t\geq t_0+T, \; \|x_t(t_0, \varphi)\| <\rho.
\end{equation}
Let
\[
S_0 = \{\varphi\in C^+: \|\varphi\|\leq \rho\}.
\]
Then, by the assumption on $f$, there is an $\alpha <0$ such that
\[
\forall i\in I_N, \forall (t, \varphi)\in {\R}_0\times S_0, \alpha \leq f_i(t, \varphi) \leq \beta_i.
\]
Define
\begin{equation}\label{e3.21}
\Omega = \{\varphi\in S_0: \forall i\in I_N, \forall \theta_1,\theta_2\in [-\tau, 0]\; \textup{with}\; \theta_1 \leq \theta_2, (\ref{e3.22})\; \textup{holds}\},
\end{equation}
\begin{equation}\label{e3.22}
\varphi_i(\theta_1) e^{\alpha(\theta_2-\theta_1)} \leq \varphi_i(\theta_2) \leq \varphi_i(\theta_1) e^{\beta_i(\theta_2-\theta_1)}.
\end{equation}
We check that this $\Omega$ is a compact uniform attractor.

(i) The compactness of $\Omega$ follows from the same proof as that of Theorem \ref{The2.1}.

(ii) For each $\varphi\in \Omega$, we have $(x_i)_t(\theta)\geq x_i(t, t_0, \varphi)e^{\beta_i\theta}$ for all $t\geq t_0$. So (\ref{e3.17}) and (\ref{e3.19}) hold for all $t_0\in {\R}_0$ and all $t\geq t_0$. Thus, $x_t(t_0, \varphi)\in S_0$ for all $t_0\in {\R}_0$ and all $t\geq t_0$. Then, with minor modification to the proof of Theorem \ref{The2.1}, we have $x_t(t_0, \varphi)\in \Omega$ for all $t_0\in {\R}_0$ and all $t\geq t_0$.

(iii) From (\ref{e3.20}) it follows that
\[
\forall \varphi\in C^+, \exists T=T(\varphi)>0 \; \textup{such that}\; \forall t_0\in{\R}_0, \forall t\geq t_0+T, \; x_t(t_0, \varphi) \in\Omega.
\]

(iv) For each $\varphi\in\Omega$, (\ref{e3.22}) implies that $\varphi_i(\theta_0)=0$ for some $i\in I_N$ and some $\theta_0\in [-\tau, 0]$ if and only if $\varphi_i(\theta)\equiv 0$ on $[-\tau, 0]$.

Therefore, $\Omega$ defined by (\ref{e3.21}) is a compact uniform attractor of (\ref{e1.2}).
\rule{2mm}{3mm}

\textbf{Proof of Theorem \ref{The2.3}}
Since $x'_i(t) \leq x_i(t)G^1_i(x^1_t)$ for $i\in I^1$ with the assumption (\ref{e2.6}) or (\ref{e2.7}), by Theorem \ref{The2.1} or \ref{The2.2} this subsystem has a compact uniform attractor $\Omega^1$. Since $F^2_i(t, \varphi^1)$ is bounded on ${\R}_0\times \Omega^1$, there are $\beta'_i>0$ for $i\in I^2$ such that $F^2_i(t, \varphi^1) \leq \beta'_i$ on ${\R}_0\times \Omega^1$. Then $x'_i(t) \leq x_i(t)(\beta'_i +G^2_i(x^2_t))$ for $i\in I^2$ and $(t_0, \varphi)\in {\R}_0\times C^+$ with $\varphi^1\in \Omega^1$. By the assumption (\ref{e2.6}) or (\ref{e2.7}) and Theorems \ref{The2.1} and \ref{The2.2}, the subsystem for $x^2$ has a compact uniform attractor $\Omega^2$. Repeating the above process, we obtain a compact uniform attractor $\Omega^k$ of the subsystem for $x^k$ for each $k\in I_m$. Then it can be verified that $\Omega^1\times \cdots \times\Omega^m$ is a compact uniform attractor for (\ref{e1.2}).
\rule{2mm}{3mm}

\section{Proof of partial permanence and permanence}

\textbf{Proof of Theorem \ref{The2.4}}
By condition (i), (\ref{e1.2}) has a compact uniform attractor $\Omega$. We first show the existence of $\rho>0$ such that
\begin{equation}\label{e4.1}
\|x_t(t_0, \varphi)-x_t(t_0, \psi)\| \leq \|\varphi-\psi\|e^{\rho(t-t_0)}
\end{equation}
for all $t_0\in{\R}_0$, $\varphi, \psi \in \Omega$ and $t\geq t_0$. From (\ref{e1.2}) we have
\begin{eqnarray*}
\qquad \qquad &\quad & x_i(t, t_0, \varphi)-x_i(t, t_0, \psi) = \varphi_i(0)-\psi_i(0) \\
    & & \qquad + \int^t_{t_0}[x_i(s, t_0, \varphi)-x_i(s, t_0, \psi)]f_i(s, x_s(t_0, \varphi))ds \\
 && \qquad + \int^t_{t_0}x_i(s, t_0, \psi)[f_i(s, x_s(t_0, \varphi))-f_i(s, x_s(t_0, \psi))]ds.
\end{eqnarray*}
By conditions (i) and (ii), there is a $\rho>0$ independent of $t_0$, $t$, $\varphi$ and $\psi$ such that
\[
|x_i(t, t_0, \varphi)-x_i(t, t_0, \psi)| \leq |\varphi_i(0)-\psi_i(0)| +\rho\int^t_{t_0}\|x_s(t_0, \varphi)-x_s(t_0, \psi)\|ds
\]
for all $t\geq t_0$ and all $i\in I_N$,  so
\[
\forall t\geq t_0, \; \|x_t(t_0, \varphi)-x_t(t_0, \psi)\| \leq \|\varphi-\psi\| +\rho\int^t_{t_0}\|x_s(t_0, \varphi)-x_s(t_0, \psi)\|ds.
\]
Then (\ref{e4.1}) follows from this and Gronwall's inequality.

By condition (iv), for each $(t_0, \varphi)\in {\R}_0\times (\cup_{i\in J}\Omega_i)$, there is a $T(t_0, \varphi)>0$ such that
    \begin{equation}\label{e4.2}
    \delta(t_0, \varphi) = \frac{1}{T(t_0, \varphi)}\int^{T(t_0, \varphi)}_0\sum_{i\in J}q_if_i(t_0+s, x_{t_0+s}(t_0, \varphi))ds >0.
    \end{equation}
Then, by (\ref{e4.1}) and the continuous dependence of $x_t(t_0, \varphi)$ on $(t_0, \varphi)$, there is an open interval $I(t_0, \varphi)\subset {\R}_0$ and an open ball $B(t_0, \varphi)$ of $\Omega$ such that
\begin{equation}\label{e4.3}
\frac{1}{T(t_0, \varphi)}\int^{T(t_0, \varphi)}_0 \sum_{i\in J}q_if_i(t_1+s, x_{t_1+s}(t_1, \psi))ds \geq \frac{1}{2}\delta(t_0, \varphi)
\end{equation}
for all $(t, \psi)\in I(t_0, \varphi)\times B(t_0, \varphi)$. Since $f(t, \varphi)$ is $T_0$-periodic by condition (iii), we may assume that, for any integer $k$ satisfying $t_0+kT_0\in {\R}_0$, $T(t_0+kT_0, \varphi) = T(t_0, \varphi)$ so $\delta(t_0+kT_0, \varphi) = \delta(t_0, \varphi)$. By (\ref{e4.1}), (ii) and (iii), we may also assume that
\begin{equation}\label{e4.4}
I(t_0+kT_0, \varphi) = I(t_0, \varphi)+kT_0, \quad B(t_0+kT_0, \varphi) = B(t_0, \varphi).
\end{equation}
Then, for any fixed $\ell\in{\R}_0$, $[\ell, \ell+T_0]\times\{\varphi\}$ is a compact set of ${\R}_0\times\Omega$ and
\[
\{I(t_0, \varphi)\times B(t_0, \varphi): t_0\in [\ell, \ell+T_0]\}
\]
is an open cover of $[\ell, \ell+T_0]\times\{\varphi\}$. Thus, there is a finite open cover of $[\ell, \ell+T_0]\times\{\varphi\}$. Combining (\ref{e4.3}) and (\ref{e4.4}) with this finite open cover, we obtain an open ball $B(\varphi)$ of $\Omega$, positive numbers $\delta_1(\varphi), \ldots, \delta_m(\varphi), T_1(\varphi), \ldots, T_m(\varphi)$, and a finite open cover $\{I^1, \ldots, I^m\}$ of ${\R}_0$ such that for each $k\in I_m$ and for all $(t_1, \psi)\in I^k\times B(\varphi)$,
\begin{equation}\label{e4.5}
\frac{1}{T_k(\varphi)}\int^{T_k(\varphi)}_0 \sum_{i\in J}q_if_i(t_1+s, x_{t_1+s}(t_1, \psi))ds \geq \delta_k(\varphi).
\end{equation}
Since $\cup_{i\in J}\Omega_i$ is compact and $\{B(\varphi): \varphi\in \cup_{i\in J}\Omega_i\}$ is an open cover of $\cup_{i\in J}\Omega_i$, there are $\varphi^1, \ldots, \varphi^p\in \cup_{i\in J}\Omega_i$ such that $\{B(\varphi^j): j\in I_p\}$ is a finite open cover of $\cup_{i\in J}\Omega_i$. Then, for each $j\in I_p$, there is an integer $m_j>0$ such that (\ref{e4.5}) holds after the replacement of $m, I^k, \varphi$ by $m_j, I^{jk}, \varphi^j$ respectively. Now put
\begin{eqnarray*}
T^0 &=& \min\{T_k(\varphi^j): j\in I_p, k\in I_{m_j}\}, \\
T^1 &=& \max\{T_k(\varphi^j): j\in I_p, k\in I_{m_j}\}, \\
\delta_0 &=& \min\{\delta_k(\varphi^j): j\in I_p, k\in I_{m_j}\}.
\end{eqnarray*}
Then the function $V: {\R}^N_+ \to {\R}_+$ defined by
\[
V(x) = \prod_{i\in J}x^{q_i}_i
\]
is continuous and $V(x) = 0$ if and only if $x\in \cup_{i\in J}\pi_i$. Thus, $x\in {\R}^N_+$ is close to $\cup_{i\in J}\pi_i$ if and only if $V(x)$ is small. By the properties of $\Omega$, $\varphi(0)$ is close to $\cup_{i\in J}\pi_i$ if and only if $\varphi\in \Omega$ is close to $\cup_{i\in J}\Omega_i$. Then we can choose $\mu>0$ sufficiently small so that the set
\[
S_1 = \{\varphi\in \Omega: 0\leq V(\varphi(0))\leq \mu\}
\]
is contained in $\cup_{j\in I_p}B(\varphi^j)$.

We claim that for each $\varphi\in S_1$ with $\varphi_i(0) >0$ for all $i\in J$ and every $t_0\in {\R}_0$, there is a $t>t_0$ such that $V(x(t, t_0, \varphi))>\mu$. Indeed, if $V(x(t, t_0, \varphi)\leq \mu$ for all $t\geq t_0$, then $x_t(t_0, \varphi)\in S_1 \subset \cup_{j\in I_p}B(\varphi^j)$ for all $t\geq t_0$. Since  $\varphi\in S_1$, we have $\varphi\in B(\varphi^j)$ for some $j\in I_p$. As $t_0\in I^{jk}$ for some $k\in I_{m_j}$, by (\ref{e4.5}) and the definition of $T^0$ and $\delta_0$, we have
\[
\int^{T_k(\varphi^j)}_0\sum_{i\in J}q_if_i(t_0+s, x_{t_0+s}(t_0, \varphi))ds \geq \delta_k(\varphi^j)T_k(\varphi^j) \geq \delta_0T^0.
\]
Differentiation of $V(x(t, t_0, \varphi))$ gives
\[
V(x(t, t_0, \varphi))' = V(x(t, t_0, \varphi))\sum_{i\in J}q_if_i(t, x_t(t_0, \varphi)).
\]
Then, with $t_1 = t_0 + T_k(\varphi^j)$, we obtain
\begin{eqnarray*}
V(t_1, t_0, \varphi)) &=& V(\varphi(0))\exp\biggl(\int^{T_k(\varphi^j)}_0\sum_{i\in J}q_if(i(t_0+s, x_{t_0+s}(t_0, \varphi))ds\biggr) \\
&\geq & V(\varphi(0))e^{\delta_0T^0}.
\end{eqnarray*}
Since $\psi=x_{t_1}(t_0, \varphi)\in B(\varphi^n)$ for some $n\in I_p$ and $t_1\in I^{nk}$ for some $k\in I_{m_n}$, by the same procedure as above and with $t_2 = t_1 + T_k(\varphi^n)$, we obtain
\[
V(x(t_2, t_0, \varphi) = V(x(t_2, t_1, \psi))\geq V(\psi(0))e^{\delta_0T^0} \geq V(\varphi(0))e^{2\delta_0T^0}.
\]
Repetition of the above process infinitely many times leads to the unboundedness of $V(x(t, t_0, \varphi))$ for $t\geq t_0$. This contradiction to our assumption $V(x(t, t_0, \varphi))\leq \mu$ for $t\geq t_0$ shows our claim.

Let $\alpha = \inf\{\sum_{i\in J}q_if_i(t, \varphi): (t, \varphi)\in {\R}_0\times \Omega\}$. The boundedness of $f$ on ${\R}_0\times\Omega$ implies $\alpha\in{\R}$. If $\alpha \geq 0$, then, for any $\varphi\in\Omega$ with $V(\varphi(0))>\mu$,
\begin{eqnarray*}
V(x(t, t_0, \varphi)) &=& V(\varphi(0))\exp\biggl(\int^t_{t_0}\sum_{i\in J}q_if_i(s, x_s(t_0, \varphi))ds\biggr) \\
    & \geq & V(\varphi(0))e^{\alpha(t-t_0)}> \mu
\end{eqnarray*}
for all $t_0\in{\R}_0$ and $t\geq t_0$.

If $\alpha < 0$ then $\rho =\mu e^{\alpha T^1}\in (0, \mu)$.
We show that $V(x(t, t_0, \varphi))>\rho$ for all $\varphi\in \Omega$ with $V(\varphi(0))>\mu$ and all $t_0\in{\R}_0$ and $t\geq t_0$. In fact, for fixed $\varphi$ and $t_0$, we have either $V(x(t, t_0, \varphi))>\mu$ for all $t\geq t_0$ or $V(x(t_1, t_0, \varphi))= \mu$ for some $t_1>t_0$ but $V(x(t, t_0, \varphi))>\mu$ for all $t\in [t_0, t_1)$. In the latter case, as $\psi = x_{t_1}(t_0, \varphi)\in S_1$, for some $j\in I_p$ and $k\in I_{m_j}$ we have
\[
\int^{T_k(\varphi^j)}_0\sum_{i\in J}q_if_i(t_1+s, x_{t_1+s}(t_1, \psi))ds >\delta_0 T^0.
\]
Thus, with $t_2 = t_1 +T_k(\varphi^j)$,
\[
V(x(t_2, t_1, \psi)) \geq V(\psi(0))e^{\delta_0T^0} >V(\psi(0))=V(x(t_1, t_0, \varphi))=\mu>\rho.
\]
For $t\in [t_1, t_2)$,
\begin{eqnarray*}
V(x(t, t_1, \psi)) &=& V(\psi(0))\exp\biggl(\int^t_{t_1}\sum_{i\in J}q_if_i(s, x_s(t_1, \psi))ds\biggr) \\
    &\geq & \mu e^{\alpha(t-t_0)}> \mu e^{\alpha T^1} = \rho.
\end{eqnarray*}
Hence, $V(x(t, t_0, \varphi))>\rho$ for all $t\geq t_0$.

This shows that for each $\varphi\in \Omega$ with $\varphi_i(0)>0$ for all $i\in J$ and for every $t_0\in {\R}_0$, there is a $T> t_0$ such that $V(x(t, t_0, \varphi))>\rho \; (\mu)$, if $\alpha <0 \; (\geq 0)$, for all $t\geq T$. Let
\begin{eqnarray}
\delta &=& \inf\{\varphi_i(0): V(\varphi(0))=\rho, \varphi\in\Omega, i\in J\}, \label{e4.6}\\
M &=&  \sup\{\|\varphi\|: \varphi\in \Omega\},\label{e4.7}
\end{eqnarray}
if $\alpha<0$ and replace $\rho$ by $\mu$ in (\ref{e4.6}) if $\alpha \geq 0$. Then
\[
\forall i\in J, \forall t\geq T, \delta \leq x_i(t, t_0, \varphi) \leq M.
\]
Therefore, (\ref{e1.2}) is partially permanent with respect to $J$.
\rule{2mm}{3mm}

\textbf{Proof of Theorem \ref{The2.5}}
By Theorem \ref{The2.4}, we need only prove that for all $(t_0, \varphi)\in {\R}_0\times (\cup_{i\in J}\Omega_i)$, there exists a $T=T(t_0, \varphi) >0$ such that
\begin{equation}\label{e4.8}
\frac{1}{T}\int^T_0\sum_{i\in J}q_i[r_i(t_0+s)+L_i(x_{t_0+s}(t_0, \varphi))]ds >0.
\end{equation}
We proceed by induction on the number $m$ of positive components of $\varphi$.

When $m=1$, we have $\varphi_k(\theta)>0$ for some $k\in I_N$ and all $\theta\in [-\tau, 0]$ and $\varphi_j(\theta)\equiv 0$ for $j\in I_N\setminus\{k\}$. Then $x_j(t, t_0, \varphi)\equiv 0$ for $t\geq t_0$ and $j\not= k$ and $x_k(t, t_0, \varphi)$ satisfies
\[
x'_k(t) = x_k(t)\biggl[r_k(t)+a_{kk}\int^0_{-\tau}x_k(t+\theta)d\xi_{kk}(\theta)- b_{kk}\int^0_{-\tau}x_k(t+\theta)d\eta_{kk}(\theta)\biggr].
\]
It can be shown (see \cite[Lemma 6]{Ho4}) that $\liminf_{t\to\infty}x_k(t) >0$. Then $\ln x_k(t)$ is bounded. Now integration of the above equation gives
\[
\frac{\ln x_k(t) -\ln x_k(t_0)}{t-t_0} = m(r_k, t_0, t)+(a_{kk}-b_{kk})m(x_k, t_0, t) + o(1)
\]
as $t\to+\infty$, where the $o(1)$ term has the precise expression
\[
\frac{1}{t-t_0}\int^0_{-\tau}\biggl[\int^{t_0}_{t_0+\theta}x_k(s)ds-\int^t_{t+\theta}x_k(s)ds\biggr]d[a_{kk}\xi_{kk}(\theta) -b_{kk}\eta_{kk}(\theta)].
\]
As the left-hand side vanishes when $t\to+\infty$ and $\lim_{t\to+\infty}m(r_k, t_0, t) = \bar{r}_k$, we must have
\[
\lim_{t\to+\infty}m(x_k, t_0, t) = \frac{\bar{r}_k}{b_{kk}-a_{kk}} > 0.
\]
Then $\lim_{t\to+\infty}m(x, t_0, t) = \hat{x}$ with $\hat{x}_k = \bar{r}_k/(b_{kk}-a_{kk})$ and $\hat{x}_j = 0$ for all $j\in I_N\setminus\{k\}$ and
\begin{eqnarray*}
\qquad & & \lim_{t\to+\infty}\frac{1}{t-t_0}\int^t_{t_0}\sum_{i\in J}q_i[r_i(s)+L_i(x_s(t_0, \varphi))]ds \\
    &=& \lim_{t\to+\infty}\sum_{i\in J}q_i[m(r_i, t_0, t)+(A-B)_im(x, t_0, t)] \\
    &=& \sum_{i\in J}q_i[\bar{r}_i +(A-B)_i\hat{x}].
\end{eqnarray*}
As $\hat{x} \in \cup_{i\in J}\pi_i$ is a fixed point of (\ref{e2.14}), by (\ref{e2.16}) we have $\sum_{i\in J}q_i(\bar{r}_i+(A-B)_i\hat{x})>0$. Then (\ref{e4.8}) holds for large enough $T>0$ when $m=1$.

Assume that (\ref{e4.8}) holds for some $m\geq 1$ and all $\varphi\in \cup_{i\in J}\Omega_i$ with at most $m$ positive components. Now suppose $\varphi^0\in \cup_{i\in J}\Omega_i$ has $m+1$ positive components and we show that (\ref{e4.8}) also holds. Let $J_1 = \{j\in I_N: \varphi^0_j(0)>0\}$ with $|J_1| = m+1$ and let
\[
\Omega^{J_1} = \cap_{j\in I_N\setminus J_1}\Omega_j.
\]
Note that $\Omega^{J_1}\subset \cup_{i\in J}\Omega_i$. Since $\Omega$ is a compact uniform attractor of (\ref{e1.2}), for each $j\in I_N$, $\Omega_j$ is a compact uniform attractor of the $N-1$-dimensional subsystem of (\ref{e1.2}) with $x_j\equiv 0$ and $\Omega^{J_1}$ is a compact uniform attractor of the corresponding $(m+1)$-dimensional subsystem. Since $\varphi^0\in \textup{int}\Omega^{J_1}$, we have $x_t(t_0, \varphi^0)\in \textup{int}\Omega^{J_1}$ for all $t\geq t_0$. There are two possible cases for the limit set $\omega(t_0, \varphi^0)$ of $x_t(t_0, \varphi^0)$ as $t\to+\infty$: (a) $\omega(t_0, \varphi^0)\not\subset \partial \Omega^{J_1}$ and (b) $\omega(t_0, \varphi^0)\subset \partial \Omega^{J_1}$.

(a) In this case, there is a $\psi\in\textup{int}\Omega^{J_1}$ and a sequence $\{t_n\}$ with $t_n\to+\infty$ as $n\to\infty$ such that $\lim_{n\to\infty}x_{t_n}(t_0, \varphi^0) = \psi$. Then $\lim_{n\to\infty}x(t_n, t_0, \varphi^0)=\psi(0)$ so the set $\{\ln x_j(t_n): j\in J_1, n\geq 1\}$ is bounded and
\begin{equation}\label{e4.10}
\forall j\in J_1, \; \lim_{n\to\infty}\frac{\ln x_j(t_n) - \ln \varphi^0_j(0)}{t_n-t_0} = 0.
\end{equation}
Integrating the $j$th component equation of (\ref{e1.2}) with (\ref{e2.9})--(\ref{e2.13}) as we did in the base case $m=1$, we obtain
\[
\frac{\ln x_j(t_n) - \ln \varphi^0_j(0)}{t_n-t_0} = m(r_j, t_0, t_n) +(A-B)_jm(x, t_0, t_n) +o(1)
\]
as $n\to\infty$. This, together with (\ref{e4.10}) and $\lim_{n\to\infty}m(r_j, t_0, t_n) = \bar{r}_j$, gives
\[
\forall j\in J_1, \; \lim_{n\to\infty}(A-B)_jm(x, t_0, t_n) = -\bar{r}_j.
\]
By choosing a subsequence of $\{t_n\}$ is necessary, we may assume that $m(x, t_0, t_n)$ tends to $\bar{x}$ as $n\to\infty$. Then $\bar{x}\in \cap_{j\in I_N\setminus J_1}\pi_j \subset \cup_{i\in J}\pi_i$ and $(A-B)_j\bar{x} = -\bar{r}_j$ for all $j\in J_1$. Thus, $\bar{x}$ is a fixed point of (\ref{e2.14}) in $\cup_{i\in J}\pi_i$ and, by (\ref{e2.16}), $\sum_{i\in J}q_i[\bar{r}_i + (A-B)_i\bar{x}] >0$. It then follows that
\begin{eqnarray*}
\qquad & & \lim_{n\to\infty}\frac{1}{t_n-t_0}\int^{t_n}_{t_0}\sum_{i\in J}q_i[r_i(s)+L_i(x_s(t_0, \varphi^0))]ds \\
    &=& \lim_{n\to\infty}\sum_{i\in J}q_i[m(r_i, t_0, t_n)+(A-B)_im(x, t_0, t_n)] \\
    &=& \sum_{i\in J}q_i[\bar{r}_i +(A-B)_i\bar{x}] >0.
\end{eqnarray*}
Therefore, for $n$ large enough, (\ref{e4.8}) holds with $T = t_n-t_0$.

(b) For each $(t_1, \varphi)\in {\R}_0\times \partial\Omega^{J_1}$, since $\varphi$ has at most $m$ positive components, by the inductive hypothesis there is a $T(t_1, \varphi)>0$ such that
\[
\int^{T(t_1, \varphi)}_0\sum_{i\in J}q_i[r_i(t_1+s)+L_i(x_{t_1+s}(t_1, \varphi))]ds >0.
\]
For this fixed $(t_1, \varphi)$, by continuous dependence there is an open interval $I(t_1, \varphi)$ of ${\R}_0$ and an open ball $B(t_1, \varphi)$ of $\Omega^{J_1}$ such that for all $(\sigma, \psi)\in I(t_1, \varphi)\times B(t_1, \varphi)$,
\[
\frac{1}{T(t_1, \varphi)}\int^{T(t_1, \varphi)}_0\sum_{i\in J}q_i[r_i(\sigma +s)+L_i(x_{\sigma +s}(\sigma, \psi))]ds > \delta(t_1, \varphi),
\]
\[
\delta(t_1, \varphi) = \frac{1}{2T(t_1, \varphi)}\int^{T(t_1, \varphi)}_0\sum_{i\in J}q_i[r_i(t_1 +s)+L_i(x_{t_1+s}(t_1, \varphi))]ds > 0.
\]
Since $r(t)$ is $T_0$-periodic and $\Omega^{J_1}$ and $\partial \Omega^{J_1}$ are compact, by the same technique as that used in the proof of Theorem \ref{The2.4} we obtain an open set $S_0$ of $\Omega^{J_1}$ with $\partial\Omega^{J_1}\subset S_0$ and numbers $\delta_0>0$, $T^0>0$ and $T^1>T^0$ satisfying for all $(\sigma, \psi)\in {\R}_0\times S_0$,
\begin{equation}\label{e4.11}
\exists T\in [T^0, T^1] \; \textup{such that}\; \frac{1}{T}\int^T_0\sum_{i\in J}q_i[r_i(\sigma +s)+L_i(x_{\sigma +s}(\sigma, \psi))]ds > \delta_0.
\end{equation}

From (\ref{e4.11}) we see that for any $(\sigma, \psi)\in {\R}_0\times S_0$, if $x_t(\sigma, \psi)\in S_0$ for all $t\geq \sigma$ then there are $T_n\geq nT^0$ such that
\begin{equation}\label{e4.12}
\int^{T_n}_0\sum_{i\in J}q_i[r_i(\sigma +s)+L_i(x_{\sigma +s}(\sigma, \psi))]ds > nT^0\delta_0 \to +\infty \quad (n\to\infty).
\end{equation}
Now that $\omega(t_0, \varphi^0)\subset \partial\Omega^{J_1} \subset S_0$, there is a $\sigma> t_0$ such that $x_t(t_0, \varphi^0)\in S_0$ for all $t\geq \sigma$. Then there is a $T>0$ for $(t_0, \varphi^0)$ such that (\ref{e4.8}) follows from (\ref{e4.12}).

By induction, (\ref{e4.8}) holds for all $(t_0, \varphi)\in {\R}_0\times (\cup_{i\in J}\Omega_i)$. Therefore, (\ref{e1.2}) with (\ref{e2.9})--(\ref{e2.13}) is partially permanent with respect to $J$.
\rule{2mm}{3mm}

\end{document}